\documentclass[vecarrow]{svmult}

\usepackage{graphicx}         
\usepackage[bottom]{footmisc} 
\usepackage{amsmath}

\def\straightletter#1{{\mathord{\rm I\mkern-3.6mu #1}}}
\def\R{\straightletter{R}}

\begin{document}

\title*{Enhancing SPH using Moving
Least-Squares and Radial Basis Functions\footnote{This work is
supported by EPSRC grants GR/S95572/01 and GR/R76615/02.}}
\titlerunning{Enhancing SPH using MLS and RBFs}
\toctitle{Enhancing SPH using Moving Least-Squares and Radial Basis
Functions}

\author{R. Brownlee\inst{1,}\thanks{Corresponding author.}
\and P. Houston\inst{2,} \and J. Levesley\inst{1}\and S.
Rosswog\inst{3}}

\institute{Department of Mathematics, University of Leicester,
Leicester, LE1 7RH, UK,
\texttt{\textrm{\{}r.brownlee,j.levesley\textrm{\}}@mcs.le.ac.uk}
\and School of Mathematical Sciences, University of Nottingham,
Nottingham NG7 2RD, UK, \texttt{paul.houston@nottingham.ac.uk} \and
School of Engineering and Science, International University Bremen,
Campus Ring 1, D-28759 Bremen, Germany,
\texttt{s.rosswog@iu-bremen.de}}

\maketitle

\begin{abstract}
In this paper we consider two sources of enhancement for the
meshfree Lagrangian particle method \emph{smoothed particle
hydrodynamics} (SPH) by improving the accuracy of the particle
approximation. Namely, we will consider \emph{shape functions}
constructed using: moving least-squares approximation (MLS); radial
basis functions (RBF). Using MLS approximation is appealing because
polynomial consistency of the particle approximation can be
enforced. RBFs further appeal as they allow one to dispense with the
\textit{smoothing-length} -- the parameter in the SPH method which
governs the number of particles within the support of the shape
function. Currently, only ad hoc methods for choosing the
smoothing-length exist. We ensure that any enhancement retains the
conservative and meshfree nature of SPH. In doing so, we derive a
new set of variationally-consistent hydrodynamic equations. Finally,
we demonstrate the performance of the new equations on the Sod shock
tube problem.
\end{abstract}

\section{Introduction}\label{rab_sec:1}
Smoothed particle hydrodynamics (SPH) is a meshfree Lagrangian
particle method primarily used for solving problems in solid and
fluid mechanics (see~\cite{rab_Monaghan05} for a recent
comprehensive review). Some of the attractive character\-istics that
SPH possesses include: the ability to handle problems with large
deformation, free surfaces and complex geometries; truly meshfree
nature (no background mesh required); exact conservation of momenta
and total energy. On the other hand, SPH suffers from several
drawbacks: an instability in tension; difficulty in enforcing
essential boundary conditions; fundamentally based on inaccurate
kernel approximation techniques. This paper addresses the last of
these deficiencies by suggesting improved particle approximation
procedures. Previous contributions in this direction (reviewed
in~\cite{rab_Belytschko96}) have focused on corrections of the
existing SPH particle approximation (or its derivatives) by
enforcing polynomial consistency. As a consequence, the
conser\-vation of relevant physical quantities by the discrete
equations is usually lost.

The outline of the paper is as follows. In the next section we
review how SPH equations for the non-dissipative motion of a fluid
can be derived. In essence this amounts to a discretization of the
Euler equations:
\begin{equation}
    \frac{\D \rho}{\D t} = -\rho \nabla \cdot v,\qquad
    \frac{\D v}{\D t} = -\frac{1}{\rho} \nabla P,\qquad
    \frac{\D e}{\D t} = -\frac{P}{\rho} \nabla \cdot v,\label{rab_euler}
\end{equation}
where $\frac{\D}{\D t}$ is the total derivative, $\rho$, $v$, $e$
and $P$ are the density, velocity, thermal energy per unit mass and
pressure, respectively. The derivation is such that important
conservation properties are satisfied by the discrete equations.
Within the same section we derive a new set of
variationally-consistent hydrodynamic equations based on improved
particle approximation. In Sect.~\ref{rab_sec:3} we construct
specific examples -- based on moving least-squares approximation and
radial basis functions -- to complete the newly derived equations.
The paper finishes with Sect.~\ref{rab_sec:4} where we demonstrate
the performance of the new methods on the Sod shock tube
problem~\cite{rab_Sod78} and make some concluding remarks.

\begin{figure}
    \centering
    \includegraphics[height=5.25cm]{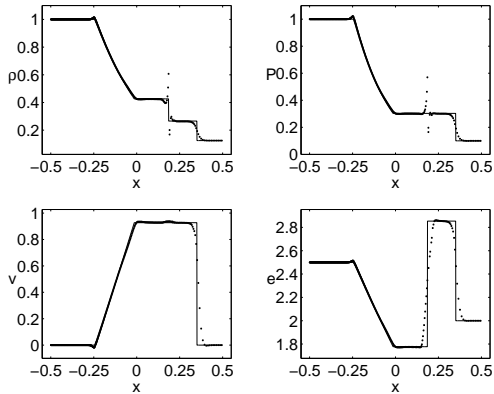}
    \caption{Shock tube simulation (t=$0.2$) using SPH.} \label{rab_fig:1}
\end{figure}

To close this section, we briefly review the SPH particle
approximation technique on which the SPH method is fundamentally
based and which we purport to be requiring improvement. From a set
of scattered particles $\{x_1,\ldots,x_N\}\subset\R^d$, SPH particle
approximation is achieved using
\begin{equation}\label{rab_sph_approx}
  Sf(x) = \sum_{j=1}^N f(x_j) \frac{m_j}{\rho_j}
  W(|x-x_j|,h),
\end{equation}
where $m_j$ and $\rho_j$ denotes the mass and density of the $j$th
particle, respectively. The function $W$ is a normalised kernel
function which approximates the $\delta$-distribution as the
\emph{smoothing-length}, $h$, tends to zero. The function
$\frac{m_j}{\rho_j} W(|x-x_j|,h)$ is called an SPH \emph{shape
function} and the most popular choice for $W$ is a compactly
supported cubic spline kernel with support $2h$. The parameter $h$
governs the extent to which contributions from neighbouring
particles are allowed to smooth the approximation to the underlying
function $f$. Allowing a spatiotemporally varying smoothing-length
increases the accuracy of an SPH simulation considerably. There are
a selection of ad hoc techniques available to accomplish this,
although often terms arising from the variation in $h$ are neglected
in the SPH method. The approximating power of the SPH particle
approximation is perceived to be poor. The SPH shape functions fail
to provide a partition of unity so that even the constant function
is not represented exactly. There is currently no approximation
theory available for SPH particle approximation when the particles
are in general positions. The result of a shock tube simulation
using the SPH equations derived in Sect.~\ref{rab_sec:2} is shown in
Fig.~\ref{rab_fig:1} (see Sect.~\ref{rab_sec:4} for the precise
details of the simulation).
The difficulty that SPH has at the contact discontinuity ($x \approx
0.2$) and the head of the rarefaction wave ($x \approx -0.25$) is
attributed to a combination of the
approximation~\eqref{rab_sph_approx} and the variable
smoothing-length not being self-consistently incorporated.

\section{Variationally-Consistent Hydrodynamic Equations}
\label{rab_sec:2}

It is well known (see \cite{rab_Monaghan05} and the references cited
therein) that the most common SPH equations for the non-dissipative
motion of a fluid can be derived using the Lagrangian for
hydrodynamics and a variational principle. In this section we review
this procedure for a particular formulation of SPH before deriving a
general set of variationally-consistent hydrodynamic equations.

The aforementioned Lagrangian is a particular functional of the
dynamical coordinates: $L(x,v) = \int \rho(v^2/2-e)\,\D x$, where
$x$ is the position, $v$ is the velocity, $\rho$ is the density, $e$
is the thermal energy per unit mass and the integral is over the
volume being discretized. Given $N$ particles $\{
x_1\ldots,x_N\}\subset \R^d$, the SPH discretization of the
Lagrangian, also denoted by $L$, is given by
\begin{equation}\label{rab_approx_lagragian}
    L = \sum_{j=1}^N m_j\Bigl(\frac{v_j^2}{2}-e_j\Bigr),
\end{equation}
where $m_j$ has replaced $\rho_j V_j$ to denote particle mass
(assumed to be constant), and $V_j$ is a volume associated with each
particle. Self-evidently, the notation $f_j$ is used to denote the
function $f$ evaluated at the $j$th particle.

The Euler-Lagrange equations give rise to SPH equations of motion
provided each quantity in~\eqref{rab_approx_lagragian} can be
written directly as a function of the particle coordinates. By
setting $f=\rho$ in~\eqref{rab_sph_approx} and evaluating at $x_j$,
we can obtain an expression for $\rho_j$ directly as a function of
the particle coordinates. Therefore, because we assume that $e_j =
e_j(\rho_j)$, the Euler--Lagrange equations are amenable.
Furthermore, in using this approach, conservation of momenta and
total energy are guaranteed via Noether's symmetry theorem. However,
when we consider improved particle approximation, the corresponding
expres\-sion for density depends on the particle coordinates in an
implicit manner, so that the Euler--Lagrange equations are not
directly amenable. To circumvent this difficulty, one can use the
principle of stationary action directly to obtain SPH equations of
motion  -- the \emph{action},
\begin{equation*}
    S=\int L\,\D t,
\end{equation*}
being the time integral of $L$. The principle of stationary action
demands that the action is invariant with respect to small changes
in the particle coordinates (i.e., $\delta S=0$). The
Euler--Lagrange equations are a consequence of this variational
principle. In~\cite{rab_Monaghan05} it is shown that if an
expression for the time rate of change of $\rho_j$ is available,
then, omitting the detail, this variational principle gives rise to
SPH equations of motion.

To obtain an expression for the time rate of change of density we
can discretize the first equation of~\eqref{rab_euler}
using~\eqref{rab_sph_approx} by collocation. By assuming that the
SPH shape functions form a partition of unity we commit error but
are able to artificially provide the discretization with invariance
to a constant shift in velocity (Galilean invariance):
\begin{equation}\label{rab_sph_den}
    \frac{\D \rho_i}{\D t} = - \rho_i \sum_{j=1}^N \frac{m_j}{\rho_j} (v_j-v_i)
    \cdot \nabla_i W(|x_i-x_j|,h_i),\qquad \text{$i=1,\ldots,N$,}
\end{equation}
where $\nabla_i$ is the gradient with respect to the coordinates of
the $i$th particle. The equations of motion that are
variationally-consistent with~\eqref{rab_sph_den} are
\begin{equation}\label{rab_sph_acc}
    \frac{\D v_i}{\D t} = -\frac{1}{\rho_i} \sum_{j=1}^N \frac{m_j}{\rho_j}\Bigl( P_i\nabla_i W(|x_i-x_j|,h_i)
    +P_j\nabla_i W(|x_i-x_j|,h_j)\Bigr),
\end{equation}
for $i=1,\ldots,N$, where $P_i$ denotes the pressure of the $i$th
particle (provided via a given equation of state). Using the first
law of thermodynamics, the equation for the rate of change of
thermal energy is given by
\begin{equation}\label{rab_thermal_energy}
  \frac{\D e_i}{\D t} = \frac{P_i}{\rho_i^2} \frac{\D \rho_i}{\D t},\qquad \text{$i=1,\ldots,N$.}
\end{equation}

As already noted, a beneficial consequence of using the
Euler--Lagrange equations is that one automatically preserves, in
the discrete equations, fundamental conservation properties of the
original system~\eqref{rab_euler}. Since we have not done this,
conservation properties are not necessarily guaranteed by our
discrete equations~\eqref{rab_sph_den}--\eqref{rab_thermal_energy}.
However, certain features of the discretization~\eqref{rab_sph_den}
give us conservation. Indeed, by virtue of~\eqref{rab_sph_den} being
Galilean invariant, one conserves linear momentum and total energy
(assuming perfect time integration). Remember that Galilean
invariance was installed under the erroneous assumption that the SPH
shape functions provide a partition of unity. Angular momentum is
also explicitly conserved by this formulation due to $W$ being
symmetric.

Now, we propose to enhance SPH by improving the particle
approxima\-tion~\eqref{rab_sph_approx}. Suppose we have constructed
shape functions $\phi_j$ that provide at least a partition of unity.
With these shape functions we form a quasi-interpolant:
\begin{equation}\label{rab_quasi_interp}
Sf = \sum_{j=1}^N f(x_j)\phi_j,
\end{equation}
which we implicitly assume provides superior approximation quality
than that provided by~\eqref{rab_sph_approx}. We defer particular
choices for $\phi_j$ until the next section. The discretization of
the continuity equation now reads
\begin{equation}\label{rab_hydro1}
    \frac{\D \rho_i}{\D t} = -\rho_i  \sum_{j=1}^N (v_j-v_i)
    \cdot \nabla \phi_j(x_i),\qquad \text{$i=1,\ldots,N$,}
\end{equation}
where, this time, we have supplied genuine Galilean invariance,
without committing an error, using the partition of unity property
of $\phi_j$. As before, the principle of stationary action provides
the equations of motion and conservation properties of the resultant
equations reflect properties present in the discrete continuity
equation~\eqref{rab_hydro1}.

To obtain~\eqref{rab_approx_lagragian}, two assumptions were made.
Firstly, the SPH shape functions were assumed to form a sufficiently
good partition of unity. Secondly, it was assumed that the kernel
approximation $\int f W(|\cdot-x_j|,h)\,\D x \approx f(x_j)$, was
valid. For our general shape functions the first of these
assumptions is manifestly true. The analogous assumption we make to
replace the second is that the error induced by the approximations
\begin{equation}\label{rab_quadrature}
  \int f \phi_j \,\D x \approx f_j \int \phi_j\,\D x \approx f_j
  V_j,\qquad \text{$j=1,\ldots,N$},
\end{equation}
is negligible. With the assumption~\eqref{rab_quadrature}, the
approximate Lagrangian associated with $\phi_j$ is identical in form
to~\eqref{rab_approx_lagragian}. Neglecting the details once again,
which can be recovered from~\cite{rab_Monaghan05}, the equations of
motion variationally-consistent with~\eqref{rab_hydro1} are
\begin{equation}\label{rab_hydro2}
 \frac{\D v_i}{\D t} = \frac{1}{m_i} \sum_{j=1}^N
     \frac{m_j}{\rho_j} P_j \nabla \phi_i(x_j),
    \qquad \text{$i=1,\ldots,N$,}
\end{equation}

The equations~\eqref{rab_thermal_energy},~\eqref{rab_hydro1}
and~\eqref{rab_hydro2} constitute a new set of
variationally-consistent hydrodynamic equations. They give rise to
the formulation of SPH derived earlier under the transformation
$\phi_j(x_i) \mapsto \frac{m_j}{\rho_j} W(|x_i-x_j|,h_i)$. The
equations of motion~\eqref{rab_hydro2} appear in~\cite{rab_Dilts00}
but along side variationally-inconsistent companion equations. The
authors advocate using a variationally-consistent set of equation
because evidence from the SPH literature
(e.g.,~\cite{rab_Bonet99,rab_Marri03}) suggests that not doing so
can lead to poor numerical results.

Linear momentum and total energy are conserved by the new equations,
and this can be verified immediately using the partition of unity
property of $\phi_j$. The $\phi_j$ will not be symmetric. However,
if it is also assumed that the shape functions reproduce linear
polynomials, namely, $\sum x_j \phi_j(x) = x$, then it is simple to
verify that angular momentum is also explicitly conserved.

\section{Moving Least-Squares and Radial Basis Functions}
\label{rab_sec:3}

In this section we construct quasi-interpolants of the
form~\eqref{rab_quasi_interp}. In doing so we furnish our newly
derived hydrodynamic
equations~\eqref{rab_thermal_energy},~\eqref{rab_hydro1}
and~\eqref{rab_hydro2} with several examples.

\subparagraph{Moving least-squares (MLS).}

The preferred construction for MLS shape func\-tions, the so-called
Backus--Gilbert approach~\cite{rab_Bos89}, seeks a quasi-interpolant
of the form~\eqref{rab_quasi_interp} such that:
\begin{itemize}
  \item $Sp = p$ for all polynomials $p$ of some fixed degree;
  \item $\phi_j(x)$, $j=1,\ldots,N$, minimise the quadratic form $
    \sum \phi_j^2(x)
    \Bigl[w\Bigl(\frac{|x-x_j|}{h}\Bigr)\Bigr]^{-1},$
\end{itemize}
where $w$ is a fixed \emph{weight function}. If $w$ is continuous,
compactly supported and positive on its support, this quadratic
minimisation problem admits a unique solution. Assuming $f$ has
sufficient smoothness, the order of convergence of the MLS
approximation~\eqref{rab_quasi_interp} directly reflects the degree
of polynomial reproduced~\cite{rab_Wendland01}.

The use of MLS approximation in an SPH context has been considered
before. Indeed, Belytschko et al.~\cite{rab_Belytschko96} have shown
that correcting the SPH particle approximation up to linear
polynomials is equivalent to an MLS approximation with
$w(|\cdot-x_j|/h)=W(|\cdot-x_j|,h)$. There is no particular reason
to base the MLS approximation on an SPH kernel. We find that MLS
approximations based on Wendland functions~\cite{rab_Wendland95},
which have half the natural support of a typical SPH kernel, produce
results which are less noisy. Dilts~\cite{rab_Dilts99,rab_Dilts00}
employs MLS approximation too. Indeed, in~\cite{rab_Dilts99}, Dilts
makes an astute observation that addresses an inconsistency that
arises due to~\eqref{rab_quadrature} -- we have the equations
\begin{equation*}
  \frac{\D V_i}{\D t} = V_i \sum_{j=1}^N  (v_j-v_i)
    \cdot \nabla \phi_j(x_i)\qquad\text{and}\qquad\frac{\D V_i}{\D t} \approx \frac{\D}{\D t} \biggl( \int \phi_i(x)\, \D
    x\biggr).
\end{equation*}
Dilts shows that if $h_i$ is evolved according to $h_i \propto
V_i^{1/d}$ then there is agreement between the right-hand sides of
these equations when a one-point quadrature of $\int \phi_i\,\D x$
is employed. Thus, providing some theoretical justification for
choosing this particular variable smoothing-length over other
possible choices.

\subparagraph{Radial basis functions (RBFs).} To construct an RBF
interpolant to an unknown function $f$ on $x_1,\ldots,x_N$, one
produces a function of the form
\begin{equation}\label{rab_rbfinterp}
    If = \sum_{j=1}^N \lambda_j \psi(|\cdot-x_j|),
\end{equation}
where the $\lambda_j$ are found by solving the linear system $
If(x_i)=f(x_i)$, $i=1,\ldots,N$. The \emph{radial basis function},
$\psi$, is a pre-specified univariate function chosen to guarantee
the solvability of this system. Depending on the choice of $\psi$, a
low degree polynomial is sometimes added to~\eqref{rab_rbfinterp} to
ensure solvability, with the additional degrees of freedom taken up
in a natural way. This is the case with the \emph{polyharmonic
splines}, which are defined, for $m>d/2$, by $\psi(|x|) = |x|^{2m-d}
\log{|x|}$ if $d$ is even and $\psi(|x|) = |x|^{2m-d}$ otherwise,
and a polynomial of degree $m-1$ is added. The choice $m\geq2$
ensures the RBF interpolant reproduces linear polynomials as
required for angular momentum to be conserved by the equations of
motion. As with MLS approximation, one has certain strong assurances
regarding the quality of the approximation induced by the RBF
interpolant (e.g. \cite{rab_Brownlee04} for the case of polyharmonic
splines).

In its present form~\eqref{rab_rbfinterp}, the RBF interpolant is
not directly amenable. One possibility is to rewrite the interpolant
in \emph{cardinal form} so that it coincides
with~\eqref{rab_quasi_interp}. This naively constitutes much greater
computational effort. However, there are several strategies for
constructing approximate cardinal RBF shape functions
(e.g.~\cite{rab_Brown05}) and fast evaluation techniques
(e.g.~\cite{rab_Beatson97}) which reduce this work significantly.
The perception of large computational effort is an attributing
factor as to why RBFs have not been considered within an SPH context
previously. Specifically for polyharmonic splines, another
possibility is to construct shape functions based on discrete
$m$-iterated Laplacians of $\psi$. This is sensible because the
continuous iterated Laplacian, when applied $\psi$, results in the
$\delta$-distribution (up to a constant). This is precisely the
approach we take in Sect.~\ref{rab_sec:4} where we employ cubic
B-spline shape functions for one of our numerical examples. The
cubic B-splines are discrete bi-Laplacians of the shifts of
$|\cdot|^3$, and they gladly reproduce linear polynomials.

In addition to superior approximation properties, using globally
supported RBF shape functions has a distinct advantage. One has
dispensed with the smoothing-length entirely. Duely, issues
regarding how to correctly vary and self-consistently incorporate
the smoothing-length vanish. Instead, a natural `support' is
generated related to the relative clustering of particles.

\section{Numerical Results}
\label{rab_sec:4}

In this section we demonstrate the performance of the
scheme~\eqref{rab_thermal_energy},~\eqref{rab_hydro1}
and~\eqref{rab_hydro2} using both MLS and RBF shape functions. The
test we have selected has become a standard one-dimensional
numerical test in compressible fluid flow -- the Sod shock
tube~\cite{rab_Sod78}. The problem consists of two regions of ideal
gas, one with a higher pressure and density than the other,
initially at rest and separated by a diaphragm. The diaphragm is
instantaneously removed and the gases allowed to flow resulting in a
rarefaction wave, contact discontinuity and shock. We set up $450$
equal mass particles in $[-0.5, 0.5]$. The gas occupying the
left-hand and right-hand sides of the domain are given initial
conditions $(P_L,\rho_L,v_L)=(1.0,1.0,0.0)$ and
$(P_R,\rho_R,v_R)=(0.1,0.125,0.0)$, respectively. The initial
condition is not smoothed.

With regards to implementation, artificial viscosity is included to
prevent the development of unphysical oscillations. The form of the
artificial viscosity mimics that of the most popular SPH artificial
viscosity and is applied with a switch which reduces the magnitude
of the viscosity by a half away from the shock. A switch is also
used to administer an artificial thermal conductivity term, also
modelled in SPH. Details of both dissipative terms and their
respective switches can be accessed through~\cite{rab_Monaghan05}.
Finally, we integrate, using a predictor--corrector method, the
equivalent hydrodynamic equations
\begin{equation}\label{rab_dVdt}
    \frac{\D V_i}{\D t} = V_i  \sum_{j=1}^N (v_j-v_i)
    \cdot \nabla \phi_j(x_i),
\end{equation}
\begin{equation*}
    \frac{\D v_i}{\D t} = \frac{1}{m_i} \sum_{j=1}^N
     V_j P_j \nabla \phi_i(x_j),\qquad
    \frac{\D e_i}{\D t} = -\frac{P_i}{m_i} \frac{\D V_i}{\D
    t},\nonumber
\end{equation*}
together with $\frac{\D x_i}{\D t} = v_i$, to move the particles. To
address the consistency issue regarding particle volume mentioned
earlier -- which is partially resolved by evolving $h$ in a
particular way when using MLS approximation -- we periodically
update the particle volume predicted by~\eqref{rab_dVdt} with $\int
\phi_i\,\D x$ if there is significant difference between these two
quantities. To be more specific, the particle volume $V_i$ is
updated if $|V_i-\int \phi_i \,\D x |/V_i \geq 1.0 \times 10^{-3}$.

\begin{figure}
\centering
\includegraphics[height=5.25cm]{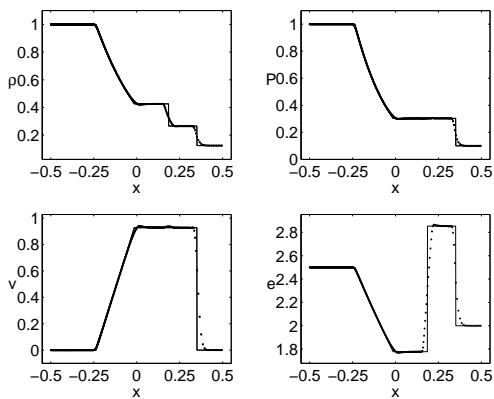}
\caption{Shock tube simulation (t=$0.2$) using linearly complete MLS
shape functions.} \label{rab_fig:2}
\end{figure}

\begin{figure}
\centering
\includegraphics[height=5.25cm]{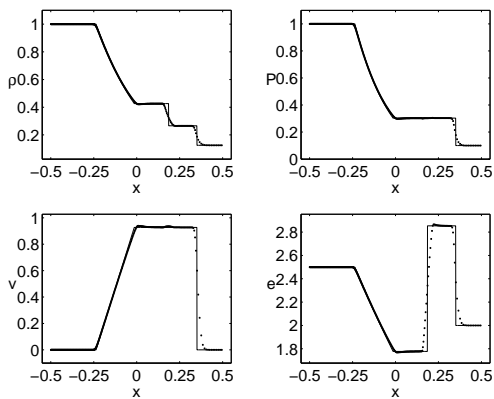}
\caption{Shock tube simulation (t=$0.2$) using cubic B-spline shape
functions.} \label{rab_fig:3}
\end{figure}

We first ran a simulation with linearly complete MLS shape
functions. The underlying univariate function, $w$, was selected to
be a Wendland function with $C^4$-smoothness. The smoothing-length
was evolved by taking a time derivative of the relationship $h_i
\propto V_i$ and integrating it alongside the other equations, the
constant of proportionality was chosen to be $2.0$. The result is
shown in Fig.~\ref{rab_fig:2}.
The agreement with the analytical solution (solid line) is
excellent, especially around the contact discontinuity and the head
of the rarefaction wave. Next, we constructed RBF shape functions.
As we mentioned in Sect.~\ref{rab_sec:3}, for this one-dimensional
problem we employ cubic B-spline because they constitute discrete
bi-Laplacians of the shifts of the globally supported basis
function, $\psi = |\cdot|^3$. The result of this simulation is shown
in Fig.~\ref{rab_fig:3}.
Again, the agreement with the analytical solution is excellent.

In the introduction an SPH simulation of the shock tube was
displayed (Fig.~\ref{rab_fig:1}). There, we
integrated~\eqref{rab_sph_den}--\eqref{rab_thermal_energy} and $h$
was updated by taking a time derivative of the relationship $h_i =
2.0m_i/\rho_i$. To keep the comparison fair, the same initial
condition, particle setup and dissipative terms were used. As
previously noted, this formulation of SPH performs poorly on this
problem, especially around the contact discontinuity. Furthermore,
we find that this formulation of SPH does not converge in the
$L_\infty$-norm for this problem. At a fixed time ($t=0.2$),
plotting number of particles, $N$, versus $L_\infty$-error in
pressure, in the region of the computational domain where the
solution is smooth reveals an approximation order of around $2/3$,
attributed to the low regularity of the analytical solution, for the
MLS and RBF methods, whereas our SPH simulation shows no
convergence. This is not to say that SPH can not perform well on
this problem. Indeed, Price~\cite{rab_Price04} shows that, for a
formulation of SPH where density is calculated via summation and
variable smoothing-length terms correctly incorporated, the
simulation does exhibit convergence in pressure. The SPH formulation
we have used is fair for comparison with the MLS and RBF methods
since they all share a common derivation. In particular, we are
integrating the continuity equation in each case.

To conclude, we have proposed a new set of discrete conservative
variation\-ally-consistent hydrodynamic equations based on a
partition of unity. These equations, when actualised with MLS and
RBF shape functions, outperform the SPH method on the shock tube
problem. Further experimentation and numerical analysis of the new
methods is a goal for future work.

\vspace{\baselineskip}
\subparagraph{\emph{\textbf{Acknowledgements.}}} The first author
would like to acknowledge Joe Mon\-aghan, whose captivating lectures
at the XIth NA Summer School held in Durham in July 2004 provided
much inspiration for this work. Similarly, the first author would
like to thank Daniel Price for his useful discussions, helpful
suggestions and hospitality received during a visit to Exeter in May
2005.

\bibliographystyle{plain}


\end{document}